\documentclass[12pt,leqno]{amsart} 
\newtheorem{pr}{Proposition}[section]
\newtheorem{rem}{Remark}[section]
\newcommand{\be}{\begin{equation}}
\newcommand{\ee}{\end{equation}}
\newcommand{\bea}{\begin{eqnarray}}
\newcommand{\eea}{\end{eqnarray}}
\newcommand{\beb}{\begin{eqnarray*}}
\newcommand{\eeb}{\end{eqnarray*}}
\usepackage{amssymb,amsfonts,amsthm,setspace,indentfirst,mathrsfs}
\usepackage{minibox}
\usepackage{enumerate}
\numberwithin{equation}{section}
\begin{document}
%
\title[A note on the existence of proper weakly cyclic $Z$ symmetric manifolds]{\bf{A note on the existence of proper weakly cyclic $Z$ symmetric manifolds}}
\author[Absos Ali Shaikh and Haradhan Kundu]{Absos Ali Shaikh and Haradhan Kundu}
\date{\today}
\address{\noindent\newline Department of Mathematics,\newline The University of 
Burdwan, Golapbag,\newline Burdwan-713104,\newline West Bengal, India}
\email{aask2003@yahoo.co.in, aashaikh@math.buruniv.ac.in}
\email{kundu.haradhan@gmail.com}
\dedicatory{}
\begin{abstract}
The object of the present note is to discuss about the defining condition of weakly cyclic  Ricci symmetric manifolds \cite{SJ06} and weakly cyclic $Z$-symmetric manifolds \cite{DMS15} and the existence of such notion by proper examples.
\end{abstract}
%
\maketitle
%
\section{Introduction}
During the study of Einstein like manifolds, by the decomposition of covariant derivative of Ricci tensor, Gray \cite{GRAY} introduced two new classes of semi-Riemannian manifold, namely, manifolds of cyclic parallel Ricci tensor (or class $\mathcal A$) and manifolds of Codazzi type Ricci tensor (or class $\mathcal B$) lying between the class of Ricci symmetric manifolds and the class of manifolds with constant scalar curvature. Again generalizing the notion of class $\mathcal A$, weakly Ricci symmetric manifolds by Tam$\acute{\mbox{a}}$ssy and Binh \cite{TQ} and generalized pseudo Ricci symmetric manifolds by Chaki \cite{Ch88}, in 2006 Shaikh and Jana \cite{SJ06} introduced the notion of weakly cyclic Ricci symmetric manifolds and studied its geometric properties. In 2015 De et al. \cite{DMS15} studied weakly cyclic $Z$-symmetric manifolds by considering the generalized $Z$-tensor of Mantica and Molinari \cite{MM12} and claimed their existence by an example. These two notions was latter studied by many authors (see, \cite{DMMS16}, \cite{Kim18}, \cite{SCJ08}, \cite{SDJ11}). 
The present note deals with the reduced form of the defining condition of weakly cyclic Ricci symmetric and weakly cyclic $Z$-symmetric manifolds along with their existence. The most important fact that the claim about the existence of weakly cyclic Ricci symmetric manifolds in \cite{SJ06} and the existence of weakly cyclic $Z$-symmetric manifolds in \cite{DMS15} are not true. The present note provides the existence of weakly cyclic Ricci symmetric and weakly cyclic $Z$-symmetric manifolds by means of three interesting examples.
%
\section{Weakly cyclic $Z$-symmetric manifolds}
Let $M$ be an $n$-dimensional ($n \ge 3$) connected smooth semi-Riemannian manifold equipped with a semi-Riemannian metric $g$. Let $\nabla$, $R$, $S$ and $r$ be respectively the Levi-Civita connection, Riemann-Christoffel curvature tensor, Ricci tensor and scalar curvature of $M$. Recently in \cite{MM12}, Mantica and Molinari defined a generalized $Z$-tensor as a symmetric (0, 2) tensor of the form $Z = S + \phi g$, where $\phi$ is an arbitrary scalar function.\\
\indent In 2015, De et al. \cite{DMS15} have studied weakly cyclic $Z$-symmetric manifolds. A Riemannian manifold $M$ is said to be weakly cyclic $Z$ symmetric \cite{DMS15} if 
\bea\label{wczs}
&&(\nabla_X Z)(X_1,X_2) + (\nabla_{X_1} Z)(X,X_2) + (\nabla_{X_2} Z)(X_1,X)\\\nonumber
&&=A(X)\, Z(X_1,X_2) + B(X_1)\, Z(X,X_2) + D(X_2)\, Z(X_1,X),
\eea
for three 1-forms $A$, $B$ and $D$ on $M$, called the associated 1-forms. Such an $n$-dimensional manifold is denoted by $(WCZS)_n$ and it is called proper if it is not weakly $Z$-symmetric \cite{MM12}. If $\phi\equiv 0$, then \eqref{wczs} reduces to
\bea\label{wcrsp}
&&(\nabla_X S)(X_1,X_2) + (\nabla_{X_1} S)(X,X_2) + (\nabla_{X_2} S)(X_1,X)\\\nonumber
&&=A(X)\, S(X_1,X_2) + B(X_1)\, S(X,X_2) + D(X_2)\, S(X_1,X),
\eea
and the manifold reduces to a weakly cyclic Ricci symmetric manifold (briefly, $(WCRS)_n$) introduced by Shaikh and Jana \cite{SJ06} in 2006. A $(WCRS)_n$ is said to be proper if it is not weakly Ricci symmetric. We mention that Shaikh and his coauthors (\cite{SCJ08}, \cite{SDJ11}) studied the weakly cyclic Ricci symmetric manifolds with some additional conditions. It may be noted that cyclic pseudo Ricci symmetric manifolds are also studied by Shaikh and Hui (\cite{SH09}, \cite{SH10a}, \cite{SH10b}).\\
\indent In the AMS Mathematical Review of \cite{SCJ08} [MR2640817], the reviewer Dacko has remarked that the defining condition \eqref{wcrsp} of a $(WCRS)_n$ reduces to 
\bea\label{wcrs}
&&(\nabla_X S)(X_1,X_2) + (\nabla_{X_1} S)(X,X_2) + (\nabla_{X_2} S)(X_1,X)\\\nonumber
&&=E(X)\, S(X_1,X_2) + E(X_1)\, S(X,X_2) + E(X_2)\, S(X_1,X),
\eea
where $E = \frac{1}{3} (A+B+D)$. Recently, Shaikh et al. \cite{SDHJK15} showed the following:\\
\textbf{Theorem A:} \cite{SDHJK15}
Let $T$ be a $(0,k)$-tensor, $k\ge 2$, on a semi-Riemannian manifold $M$, symmetric with respect to $p$-number ($2p\le k$) of indices with other $p$-number of indices taken together (including non-identical arrangement). If the manifold satisfies
\beb
(\nabla_X T)(X_1,X_2,\cdots,X_k) &=& \pi(X)T(X_1,X_2,\cdots,X_k)\\
																	& & +\pi_1(X_1)T(X,X_2,\cdots,X_k)+\cdots+\pi_k(X_k)T(X_1,X_2,\cdots,X),
\eeb
then each corresponding 1-form $\pi_i$ among those $p$ indices can be replaced by the same in pair.\\
\indent Since $(\nabla_X Z)(X_1,X_2) + (\nabla_{X_1} Z)(X,X_2) + (\nabla_{X_2} Z)(X_1,X)$ is symmetric in $X_1, X_2$ and $X_3$, motivating from the AMS review of \cite{SCJ08} and in view of Theorem A, we get the following result:
\begin{pr}\label{pr2.1}
A semi-Riemannian manifold $(M, g)$ satisfying \eqref{wczs} also satisfies
\bea\label{rf}
&&(\nabla_X Z)(X_1,X_2) + (\nabla_{X_1} Z)(X,X_2) + (\nabla_{X_2} Z)(X_1,X)\\\nonumber
&&=E(X)\, Z(X_1,X_2) + E(X_1)\, Z(X,X_2) + E(X_2)\, Z(X_1,X),
\eea
where $E = \frac{A+B+D}{3}$.
\end{pr}
\begin{rem} 
We note that there may exist some $(WCZS)_n$ (resp., $(WCRS)_n$) whose associated 1-forms $A$, $B$ and $D$ are different, i.e., the associated 1-forms of the structure $(WCZS)_n$ (resp., $(WCRS)_n$) may not unique. We also note that the metric given in Example 2 (resp., Example 1) is a $(WCZS)_4$ (resp., $(WCRS)_4$) with unique associated 1-forms. Still there is a natural problem to find out a proper $(WCZS)_n$ (resp., $(WCRS)_n$) with non-unique associated 1-forms, i.e., a $(WCZS)_n$ (resp., $(WCRS)_n$) with $A \ne B \ne D$. Although we present an improper example of $(WCRS)_n$ with $A \ne B \ne D$.
\end{rem}
\begin{rem} 
Very recently, Kim \cite{Kim18} studied some geometric properties of a $(WCZS)_n$. All the results of this paper are established by considering $B-D\ne 0$ and hence these results are standing on the problem of the existence of $(WCZS)_n$ with $A \ne B \ne D$.
\end{rem}
\begin{rem} 
The Theorem 8.1 of  \cite{DMS15} is based on the  sufficient condition for a quasi Einstein manifold to be Ricci semisymmetric. Thus it is an easy consequence of Theorem 1 of \cite{DSS06}. Again the Theorem 8.2 of \cite{DMS15} is an easy consequence of Lemma 5.6 of \cite{SK14} and Theorem 8.1.
\end{rem}
\section{Proper examples of $(WCRS)_n$ and $(WCZS)_n$}\label{exam}
The metric in Example 9.1 of \cite{DMS15} (hence Example 1 of \cite{Kim18}) is the metric given in Example 3 of \cite{SJ06}. We note that this is not a $(WCZS)_4$ for $\phi = \frac{1}{(x^4)^2}$ due to unaccented of $Z_{14,1} = -\frac{8}{9 (x^4)^{5/3}}$. If we take $\phi = \frac{2}{3(x^4)^2}$, then the Example 9.1 of \cite{DMS15} is a $(WZS)_4$ and hence improperly a $(WCZS)_4$.\\
\indent It may be mentioned that the examples given in \cite{SJ06} also do not satisfy the defining condition of weakly cyclic Ricci symmetric manifold. We now present two proper examples of  $(WCRS)_4$ and $(WCZS)_4$.\\
\textbf{Example 1:} Consider the G\"odel type metric (\cite{DHJKS14}, \cite{SK16srs}) in cylindrical coordinate system
\be
ds^{2} = (dt + h(r)\, d\phi )^{2} - (f(r))^{2}\, d\phi ^{2} - dr^{2} - dz^{2},
\ee
where $h(r) = \frac{r^2}{8}$ and $f(r) = \frac{r}{8} \sqrt{8+r^2}$.
The non-zero (upto symmetry) components of $R$ and $S$ are given by
$$R_{1212}= - \frac{r^2}{64}, \ \ \ R_{1313}= R_{1323}= R_{2323}= -\frac{1}{r^2+8},$$
$$S_{11}=-\frac{2}{r^2+8}, \ \ \ S_{12}= S_{22}= -\frac{r^2+16}{8 \left(r^2+8\right)}, \ \ \ S_{33}=\frac{8}{\left(r^2+8\right)^2}.$$
Again the non-zero (upto symmetry) components of $\nabla R$ and $\nabla S$ are given by
$$\frac{1}{2}R_{1212,3}= R_{1213,2}= -R_{1223,1}=\frac{r^3}{64 \left(r^2+8\right)}, \ \ R_{1313,3}= R_{1323,3}= R_{2323,3}=\frac{4 r}{\left(r^2+8\right)^2},$$
$$S_{11,3}= 6 S_{13,1}=\frac{6 r}{\left(r^2+8\right)^2},\ \ \ S_{12,3}= S_{22,3}=\frac{r \left(r^2+24\right)}{4 \left(r^2+8\right)^2},$$
$$S_{13,2}= S_{23,1}= S_{23,2}=\frac{r}{8 \left(r^2+8\right)}, \ \ \ S_{33,3}=-\frac{32 r}{\left(r^2+8\right)^3}.$$
It is easy to check that if we consider $A = B = D = \left\{0, 0, -\frac{4 r}{8 + r^2}, 0\right\}$, then \eqref{wcrsp} holds and hence the manifold is a $(WCRS)_4$.\\
\textbf{Example 2:} Let us consider the Lorentzian manifold $M = \{(x^1,x^2,x^3,x^4)\in \mathbb R^4:x^1>0\}$ endowed with the metric
\be
ds^{2} = (x^1)^2 (dx^1)^{2} + (x^1)^2 (dx^2)^{2} + (x^1)^2 (dx^3)^{2} - (x^1)^{-2} (dx^4)^{2}.
\ee
Then the non-zero (upto symmetry) components of $R$, $S$, $\nabla R$ and $\nabla S$ are given by
$$R_{1212}= R_{1313}= -R_{2323}=1, \ \ \frac{1}{3}R_{1414}= - R_{2424}= - R_{3434}=\frac{1}{(x^1)^4},$$
$$S_{11}= -S_{22}= -S_{33}=\frac{1}{(x^1)^2}, \ \ S_{44}=-\frac{1}{(x^1)^6},$$
$$-R_{1212,1}= -2R_{1223,3}= -R_{1313,1}= 2R_{1323,2}= R_{2323,1}=\frac{4}{x^1},$$
$$-\frac{1}{3}R_{1414,1}= R_{1424,2}= R_{1434,3}= R_{2424,1}= R_{3434,1}=\frac{4}{(x^1)^5},$$
$$-S_{11,1}= 2S_{12,2}= 2S_{13,3}= S_{22,1}= S_{33,1} = \frac{4}{(x^1)^3}, \ \ S_{44,1} = \frac{4}{(x^1)^7}.$$
Let us now consider $\phi = \frac{a}{(x^1)^6}-\frac{1}{(x^1)^4}$, $a$ is a constant, and then the non-zero (upto symmetry) components of $Z=S + \phi g$ and $\nabla Z$ are given by
$$Z_{11}=\frac{a}{(x^1)^4}, \ \ Z_{22}= Z_{33}= -\frac{2 (x^1)^2-a}{(x^1)^4}, \ \ Z_{44}= -\frac{a}{(x^1)^8},$$
$$Z_{11,1}=-\frac{6 a}{(x^1)^5}, \ \ Z_{12,2}= Z_{13,3}=\frac{2}{(x^1)^3}, \ \ Z_{22,1}= Z_{33,1}=\frac{2 \left(4 (x^1)^2-3 a\right)}{(x^1)^5}, \ \ Z_{44,1} = \frac{6 a}{(x^1)^9}.$$
If we consider $A = B = D = \left\{-\frac{6}{x^1}, 0, 0, 0\right\}$, then \eqref{wczs} holds and hence the manifold is a $(WCZS)_4$.
\begin{rem}
Very recently De et al. \cite{DMMS16} studied $(WCZS)_4$ spacetimes (connected four dimensional semi-Riemannian manifold with Lorentzian metric). They claimed that such spacetimes are quasi-Einstein (Proposition 2.1 \cite{DMMS16}). In our Example 2, $M$ is a $(WCZS)_4$ spacetime but it is not quasi-Einstein. To prove Proposition 2.1 of \cite{DMMS16}, the authors considered the vector field corresponding to $F=B-D$ as unit timelike but it is not true for any $(WCZS)_4$ spacetime, e.g., $F\equiv 0$ in the above example. Hence from Example 2 we can conclude that the Proposition 2.1 of \cite{DMMS16} is not always true.
\end{rem}
\noindent\textbf{Note:} In Example 1 (resp., Example 2), the metric is proper $(WCRS)_4$ (resp., $(WCZS)_4$) but the associated 1-forms $A,B$ and $D$ are all equal and hence unique. We now present an example of improper $(WCRS)_4$ with non-unique associated 1-forms.\\
\textbf{Example 3:} Let $M$ be a non-empty open connected subset of $\mathbb R^4$ endowed with the semi-Riemannian metric
\be
ds^{2} = e^{x^1} \left[(x^3)^2+x^3+1\right] (dx^1)^{2} + 2 dx^1 dx^2 + (dx^3)^{2} + (x^1)^2 (dx^4)^{2}.
\ee
Then the non-zero (upto symmetry) components of $R$, $S$, $\nabla R$ and $\nabla S$ are given by
$$R_{1313}= -e^{x^1}, \ \ \ R_{1313,1}= -e^{x^1}, \ \ \ \ S_{11}= e^{x^1}, \ \ \ S_{11,1} = e^{x^1}.$$
Obviously the metric is recurrent as well as Ricci recurrent and hence it is improperly a $(WCRS)_4$. One can easily check that it satisfies the this metric satisfies the weakly cyclic Ricci symmetry condition \eqref{wcrsp} for $A=\{\alpha,0,0,0\}$, $B=\{\beta,0,0,0\}$ and $D=\{3-\alpha-\beta,0,0,0\}$, where $\alpha$ and $\beta$ are any non-vanishing scalars.
%


\end{document}